\documentstyle{article}
\begin{document}
\bigskip

\begin{center}

{\Large \bf
 Computer Algebra of Vector Bundles, Foliations and Zeta Functions
and a Context of Noncommutative Geometry  
 }
\end{center}
\smallskip
\begin{center}
{\bf Nikolaj M. Glazunov} \end{center}
\smallskip
\begin{center}
{\rm Glushkov Institute of Cybernetics NAS \\
 03680 Ukraine Kiev-680 Glushkov prospekt 40 \\
 Email:} {\it glanm@d105.icyb.kiev.ua }
\end{center} \smallskip

 \begin{center} {\bf  Abstract} \end{center}
  We present some methods and results in the application of
algebraic geometry and computer algebra to the study of
algebraic vector bundles, foliations and zeta functions.
A connection of the methods and results with noncommutative
geometry will be consider.

\section*{Introduction}
 This is a work on perspective. \\
 In the survey I want to discuss algebraic and computer algebra 
aspects of vector bundles, foliations and zeta functions from
commutative and noncommutative geometry points of view.
I shall do it on the base of some parts of papers
of A. Connes~\cite{Co:N2,Co:NG}, H. Moriyoshi~\cite{Mo:OA} and 
N. Nekrasov and A. Schwarz~\cite{NS:I} on noncommutative torus, 
operator algebras, and instantons on noncommutative ${\bf R}^4$,
of A. Connes and D. Kreimer~\cite{CoK} on renormalization and
the Hopf algebra structure of graphs and D. Broadhurst and 
D. Kreimer~\cite{BK:R} on renormalization automated by Hopf
algebra. Algebraic varieties appears in many mathematical and
physical problems (strings as algebraic curves, some Calabi-Yau 
manifolds). So I want mention some novel papers on computer
algebra~\cite{Sc:I,Sc:R,DJGP} which have strong connections with  
algebraic geometric.   \\
  In section 1 we give formulae for computation of local 
charts, tangent bundles and transition functions of two dimensional 
sphere. These formulae (and some formulae
of the next subsection) can be implemented on computer
algebra systems (Reduce, Maple) straightforward. Most formulae
are tested by author on different versions of Reduce and Maple. 
Then we discuss vector bundles over algebraic curves, Kronecker
foliation, formal groups, moduli spaces and connections.   \\
  As have noted by A. Connes~\cite{Co:N2}, there is "a spectral
interpretation of zeros of zeta and L-functions in terms of
constructions involving adeles, more specifically the
noncommutative space of adele classes". In section 2 we give
formulae for computation of values of some zeta and L-functions.
 \\
  Section 3 presents very short discussion of the algebraic context
of noncommutative geometry. \\   
  For some problems it is necessary to programming them from
the beginning (elementary data stuctures, efficient algorithm). 
A useful standpoint of solving a variety of problems is to model them
in terms of some graphs (still 1959 F. Harary has proposed a
graph method for complete reduction of a matrix with a view toward
finding its eigenvalues~\cite{Ha:G}). These include trees, branching and 
connectivity with cutset of graphs, covering problems, networks
and flows, matching and maximal matching. Graphs can modeling
commutative diagrams and complexes as well as Feynman diagrams. 
Under interaction of strings manifolds of varying dimension can appear. 
 Solving differential equations on the manifolds by an iterative 
method (for instance, by Newton map $ N_{f}(x) = x - [Df(x)]^{-1}f(x)$) 
we have to compute Jacobian. The efficien algorithm can be implemented 
by the cutset method. 
In the Appendix we give our implementation of the cutset algorithm. \\
 Side by side with works of above-mentioned authors, 
the survey based on my talks on SCAN2000~\cite{G:VN},
Banach Center (Warsaw) conference~\cite{G:F}, SAGP'99 (Luminy, France),
IMACS-ACA~\cite{G:B}, Mittag-Leffler seminar "Geometry and 
Physics" (1998), Math. Institute of Stockholm University seminar 
"Algebra and Geometry"~\cite{G:M}, Dubna'98~\cite{G:C} and
on SNADE'97~\cite{G:A}.  \\
 So the purpose of the talk is two-fold: first to review
algebraic geometric and computer algebra aspects of vector 
bundles, foliations and zeta functions; second to disscuss  some 
noncommutative context of the structures and methods
in the frame of algebraic and analytic computations.

\section{Geometry, Topology and Dynamics }
As elements of differential geometry~\cite{Go} and
differential topology~\cite{Hi,Ch} play an important role in our
consideration we recall for spesialists in computer science 
notions of local charts, atlases,  differentiable manifold and
tangent bundles by the example of two dimensional sphere
$ S^{2}: x_{1}^2 + x_{2}^2 + x_{3}^2 = 1 $.  \\
{\it Local charts as half-spheres} of $ S^2 $.   \\
In the case local charts are $ (U_{i},\varphi_{i}) $, where
$$ x = (x_{1},x_{2},x_{3}), \;
 U_{i} = (x \in S^2 : x_{i} > 0), \;
 U_{3+i} = (x \in S^2 : x_{i} < 0), $$
$i = 1,2,3; \; {\bigcup}_{j=1}^{6} U_{j} = S^{2} $ and maps
$\varphi_{i}: U_{i} \rightarrow {\bf R}^{2} $ are follows \\
$$\varphi_{1}(x_{2},x_{3}): u_{1}=x_{2}, u_{2}=x_{3}; $$
$$\varphi_{2}(x_{1},x_{3}): v_{1}=x_{1}, v_{2}=x_{3}; $$
$$\varphi_{2}^{-1}(v_{1},v_{2}) = (v_{1},\sqrt{1-v_{1}^{2} -
 v_{2}^{2}}, \; v_{2}) \; ; $$
$$\varphi_{3}(x_{1},x_{2}): w_{1}=x_{1}, w_{2}=x_{2}; $$
 Maps $ \varphi_{4}, \varphi_{5}, \varphi_{6} $ are defined
by symmetry. \\
   The intersection $ U_{1} \cap U_{2} \neq \emptyset .$ So
$$ \varphi_{1}(x_{2},x_{3}) = (u_{1}, u_{2}), \;
  \varphi_{2}(x_{1},x_{3}) = (v_{1}, v_{2}). $$
and the map
$$\varphi_{1} \circ \varphi_{2}^{-1} =
(u_{2} = v_{2}, u_{1} = \sqrt{1 - v_{1}^{2} - v_{2}^{2}} \;) $$
is a $ C^{r}- $ diffeomorphism between open sets
$ \varphi_{1}(U_{1} \cap U_{2}) $ and
$ \varphi_{2}(U_{1} \cap U_{2}).$ The diffeomorphism
$\varphi_{1} \circ \varphi_{2}^{-1} $ is called the changes of
coordinates. \\
{\it Local charts and stereographics projection} of $ S^2 $.   \\
 In the case there are two charts
$$ U_{1} = (S^{2} \setminus (0,0,1) \mid
   \varphi_{1} = (\frac{x_{1}}{1 - x_{3}},
    \frac{x_{2}}{1 - x_{3}})) $$
   and
$$ U_{2} = (S^{2} \setminus (0,0,-1) \mid
   \varphi_{2} = (\frac{x_{1}}{1 + x_{3}},
    \frac{x_{2}}{1 + x_{3}})) $$
 Now
    $$ U_{1} \cap U_{2} = (x \in S^{2} \mid |x_{3}| < 1) .$$
Let $ (u_{1},u_{2}), (v_{1},v_{2}) $ be local coordinates at
$ U_{1}, U_{2} \; $ respectively.
  In the intersection  $ U_{1} \cap U_{2} $
$$  \varphi_{1}^{-1} = (\frac{2u_{1}}{u_{1}^{2} +
 u_{2}^{2} + 1},
\frac{2u_{2}}{u_{1}^{2} + u_{2}^{2} + 1},
\frac{u_{1}^{2} + u_{2}^{2} - 1}{u_{1}^{2} + u_{2}^{2} + 1}), $$
$$  \varphi_{2}^{-1} = (\frac{2v_{1}}{v_{1}^{2} +
 v_{2}^{2} + 1},
\frac{2v_{2}}{v_{1}^{2} + v_{2}^{2} + 1},
\frac{1 - v_{1}^{2} - u_{2}^{2}}{v_{1}^{2} + v_{2}^{2} + 1}), $$
and in the case \\
$\varphi_{1} \circ \varphi_{2}^{-1}(v_{1},v_{2}) = (v_{1},v_{2}),
$ and
$\varphi_{2} \circ \varphi_{1}^{-1}(u_{1},u_{2}) = (u_{1},u_{2}),
$ \\
Again we have diffeomorphisms between
$ \varphi_{1}(U_{1} \cap U_{2}) $ and
$ \varphi_{2}(U_{1} \cap U_{2}).$ \\
{\it Tangent bundles} of one and two-spheres $ S^{1} $ and
$ S^{2} \; $ . \\
 Here we recall tangent bundles of  $ S^{1} $ and
$ S^{2} \; $ . In both cases the tangent bundle is a
differentiable manifold. Tangent bundle $ TS^{1} $ of
$ S^{1} $ is the product $ S^{1} \times {\bf R} $.
Therefore $ TS^{1} $ is a trivial bundle and manifold
$ S^{1} $ is {\it parallelizable}. \\
  Tangent bundle $ TS^{2} $ of $ S^{2} $ can be defined
by 3 local charts (in the sense of vector bundles; see
below). Let $ U_{i} = (x \in S^{2} | |x_{i}| < 1). $
Then $ \; {\bigcup}_{i=1}^{3} U_{i} = S^{2}. $ Let $ E $
be the set $ (x,t) $ in $ {\bf R}^3 \times {\bf R}^3 $
with $ x \in S^{2} $ as above and $ t = (u,v,w) $
such that scalar product $ <x,t> = 0. $ Let
$ p: (x,t) \mapsto x $ be the projection $ E $ on
$ S^{2}. $ Homeomorphisms \\
$$ \Phi_{1}: \; p^{-1}(U_{1}) \to (U_{1} \times {\bf R}^{2}) $$
$$ (x_{1},x_{2},x_{3};u,v,w) \mapsto (x_{1},x_{2},x_{3};vx_{3}
- wx_{2},u), $$
$$ \Phi_{2}: \; p^{-1}(U_{2}) \to (U_{2} \times {\bf R}^{2}) $$
$$ (x_{1},x_{2},x_{3};u,v,w) \mapsto (x_{1},x_{2},x_{3};wx_{1}
- ux_{3},v), $$
$$ \Phi_{3}: \; p^{-1}(U_{3}) \to (U_{3} \times {\bf R}^{2}) $$
$$ (x_{1},x_{2},x_{3};u,v,w) \mapsto (x_{1},x_{2},x_{3};ux_{3}
- vx_{2},w), $$
define atlas $ (U_{i},\Phi_{i}) $ of tangent bundle $ TS^{2}. $
There is well known fact that $ TS^{2} $ is not trivial
bundle and therefore the manifold $ S^{2} $ is not
parallelizable.      \\
   The above tangent bundles of spheres are really vector
bundles. This means that in the case of $ S^{2} $ there are
unique continuous mappings \\
$$ g_{ij} : U_{i} \cap U_{j} \rightarrow GL(2,{\bf R}), $$
such that the mapping
$$ \psi_{ij} = \Phi_{i} \circ \Phi_{j}^{-1} :
 (U_{i} \cap U_{j}) \times {\bf R}^{2} \rightarrow
 (U_{i} \cap U_{j}) \times {\bf R}^{2}, $$
satisfies $ \psi_{ij} = (x,g_{ij}(x)y) $ for every
$ (x,y) \in (U_{i} \cap U_{j}) \times {\bf R}^{2}. $
In the case $ TS^{2} $ the {\it transition functions}
 $ g_{ij} $ are follows:
$$  g_{21}(x) = \frac{-1}{x_{2}^{2} + x_{3}^{2}} \;
  \left| \begin{array}{cc}
  x_{1}x_{2}& x_{3} \\
  -x_{3}&   x_{1}x_{2}
  \end{array} \right| \; ,
$$
$$  g_{32}(x) = \frac{-1}{x_{1}^{2} + x_{3}^{2}} \;
  \left| \begin{array}{cc}
  x_{2}x_{3}& x_{1} \\
  -x_{1}&   x_{2}x_{3}
  \end{array} \right|  \; ,
$$
$$ g_{13}(x) = \frac{-1}{x_{1}^{2} + x_{2}^{2}} \;
  \left| \begin{array}{cc}
  x_{1}x_{3}& x_{2} \\
  -x_{2}&   x_{1}x_{3}
  \end{array} \right| \; .
  $$

\subsection{Vector Bundles }
   More generally  we can give a manifold by charts and local
diffeomorphisms. 
{\it Local chart} or a {\it system of coordinates} on a
topological space $ M $ is a pair $ (U,\varphi) $ where $ U $
is an open set in $ M $ and $\varphi : \ U \rightarrow {\bf R}^{m} $
is a homeomorphism from $ U $ to an open set $\varphi(U) $ in
$ {\bf R}^{m} $. An {\it atlas} $\Phi $ of dimension $ m $ is
a collection of local charts whose domains cover $ M $ and such
that if $ (U,\varphi ), (U_{1}, \varphi_{1}) \in \Phi $ and
$ U \cap U_{1} \neq 0 $ then the map
$$ \varphi_{1} \circ \varphi^{-1}:\varphi(U \cap U_{1}) \rightarrow
 \varphi_{1} (U \cap U_{1}) $$
is a $ C^{r} $-diffeomorphism between open sets in $ {\bf R}^m $. \\
 {\it Fibre space} is the object $(E,p,B)$, where $ p $
is the continuous surjective (= on) mapping of a topological
space $ E $ onto a space $ B $ (in our consideration $ B = M $
is a differential manifold), and $ p^{-1}(b) $ is called the
{\it fibre} above $ b \in B $. Both the notation
$ p:E \rightarrow B $
and $ (E,p,B) $ are used to denote a {\it fibration}, a
{\it fibre space}, a {\it fibre bundle} or a {\it bundle}. \\
{\it Vector bundle} is fibre space each fibre $ p^{-1}(b) $ of
which is endowed with the structure of a (finite dimension)
vector space $ V $ over skew-field $ K $ such that the following
local triviality condition is satisfied: each point $ b \in B $
has an open neighborhood $ U $ and a $ V $ -isomorphism of
fibre bundles $ \phi: \ p^{-1}(U) \rightarrow U \times V $ such
that $ \phi \mid_{p^{-1}(b)}: \ p^{-1}(b) \rightarrow b \times V $
is an isomorphism of vector spaces for each $ b \in B $. \\
 $ dim V $ is said to be the dimension of the vector bundle. \\
An {\em Hermitian bundle} over algebraic variety $X$
consists of a vector bundle over $X$ and a choice of $ C^\infty$
Hermitian metric on the vector bundle over complex manifold
$ X({\bf C}) $, which is invariant under antiholomorphic involution
of $ X({\bf C}) $.            \\

{\it Tangent bundle}~\cite{Hi,Go,Ch}. \\
The {\it tangent space} to a differentiable manifold $ M $ at
point $ a \in M $ can be defined as the set of tangency classes
of smooth paths in $ M $ based at $a$. It will be denoted by
$ T_{a}M .$ Elements of $ T_{a}M $ are called tangent vectors to
$ M $ at $a$. \\
The {\it tangent bundle} of $ M$, denoted by $ TM $, is the union
of the tangent spaces at all the points of $ M .$ By well known
way $ TM $ can be made into a smooth manifold. \\
Recall  well known facts about $ TM: $ \\
(i) if $M$ is $C^{r}$ then $TM$ is $C^{r-1};$     \\
(ii) if $ M $ is $ C^{\infty} $ or $ C^{\omega} $  then the
same holds for $ TM $; \\
(iii) if $M$ has dimension $n$ then $TM$ has dimension $2n$; \\
(iv) there is a natural map $p: TM \rightarrow M $ called the
{\it projection} map, taking $ T_{a}M $ to $a$ for each $a$ in
$M$, i.e. $p$ takes all tangent vectors at $a$ to the point
$a$ itself. Thus $p^{-1}(a) = T_{a}M $ (fibre of the bundle
over $a$). The projection $p$ is a smooth map
$C^{r-1}$ if $M$ is $C^{r}; $ \\
{\bf Vector fields and Flows.} \\
A {\it vector field} on a smooth manifold $M$ is a map
$F: M \rightarrow TM $ which satisfies $ p \circ F = id_{M} $,
where $p$ is the natural projection  $ TM \rightarrow M. $
By its definition a vector field is a {\it section} of the
bundle $ TM .$  \\
Let $F$ be a vector field on $M$. A {\it solution curve}
to $F$, based at $a$ on $M$ is a path $ c: I \rightarrow M $
(where $ I $ is some open interval $ (a,b) $ with $ a < 0 < b $ )
satisfying $ c(0) = a $ and $ T_{t}c = F(c(t)) $ for all
$t$ in $I$. (Here $ TI = I \times {\bf R} $ and so $ T_{t}c $
is a map $ {\bf R} \rightarrow TM. $ This demonstrate
the well known fact that a vector field on a manifold is the
global version of a first order autonomous system of $n$
ordinary differential equations on $ {\bf R}^{n}.$ \\
Let
$$ \dot {\bf x} =  A{\bf x},   \quad       (2)   $$
where ${\bf x} \in {\bf R}^{n}$, $ A$ is $n \times n$
matrix, be a linear system of ordinary differential equations.
It is well known that the solution of the system (2) together
with the initial condition ${\bf x}(0) = {\bf x}_{0} $ is
given by
$$ {\bf x}(t) = e^{At}{\bf x}_{0}, $$
where $ e^{At} $ is an $n \times n$ matrix function defined by
its Taylor series. The mapping $ e^{At}: {\bf R}^{n}
\rightarrow {\bf R}^{n} $ is called {\it the flow of linear
system} (2). \\
Recall the definition of the flow ${\phi}_{t}$ of the
nonlinear system
$$ \dot {\bf x} =  F({\bf x}),   \quad       (3)   $$
Let $I({\bf x}_{0})$ be the maximal interval of existence
of the solution ${\phi}(t,{\bf x}_{0}) $  of (2) with
the initial value ${\bf x}(0) = {\bf x}_{0}. $ Let X
be an open subset of ${\bf R}^{n}$ and let $ F \in C^{1}(X).$
For $ {\bf x}_{0} \in X, $ let ${\phi}(t,{\bf x}_{0}) $ be
the solution of (3) with initial value problem
${\bf x}(0) = {\bf x}_{0} $ defined on its maximal
interval of existence $I({\bf x}_{0}).$ Then for $ t \in
I({\bf x}_{0}),$ the mapping ${\phi}_{t}: X \rightarrow X $
defined by
$${\phi}_{t}({\bf x}_{0}) = {\phi}(t,{\bf x}_{0}) $$
is called {\it the flow of the differential equation} (1)
or {\it the flow of the vector field} F({\bf x}).

\subsection{Vector Bundles over Projective Algebraic Curves }

Let $X$ be a projective algebraic curve over algebraically
closed field $k$ and $g$ the genus of $X$. 
Let ${\cal VB}(X)$ be the category of vector bundles over 
$X$. Grothendieck have shown that for a rational curve every
vector bundle is a direct sum of line bundles. Atiyah have classified 
vector bundles over elliptic curves. The main result is  \\
{\bf Theorem}~\cite{At}. Let $X$ be an elliptic curve, $A$ a
fixed base point on $X$. We may regard $X$ as an abelian
variety with $A$ as the zero element.  Let ${\cal E}(r,d)$ denote 
the the set of equivalence classes of indecomposable vector
bundles over $X$ of dimension $r$ and degree $d$. Then each
${\cal E}(r,d)$ may be identified with $X$ in such a way that \\
$ det: {\cal E}(r,d) \rightarrow {\cal E}(1,d) $ corresponds to 
$ H: X \rightarrow X, $ \\
where $ H(x) = hx = x + x + \cdots + x \; (h \:$ times$)$, and
$h = (r,d)$ is the highest common factor of $r$ and $d$. \\
Curve $X$ is called a {\it configuration} if its normalization
is a union of projective lines and all singular points of $X$
are simple nodes. For each configuration $X$ can assign a
non-oriented graph $\Delta(X)$, whose vertices are irreducible
components of $X$, edges are its singular and an edge is
incident to a vertex if the corresponding component contains the
singular point. Drozd and Greuel have proved: \\
{\bf Theorem}~\cite{DG}. 1. ${\cal VB}(X)$ contains finitely many
indecomposable objects up to shift and isomorphism if and only if
$X$ is a configuration and the graph $\Delta(X)$ is a simple chain
(possibly one point if $X = {\bf P}^1$). \\
2. ${\cal VB}(X)$ is tame, i.e. there exist at most one-parameter
families of indecomposable vector bundles over $X$, if and only if
either $X$ is a smooth elliptic curve or it is a configuration
and the graph $\Delta(X)$ is a simple cycle (possibly, one loop
if $X$ is a rational curve with only one simple node).   \\
3. Otherwise ${\cal VB}(X)$ is wild, i.e. for each finitely generated
$k-$algebra $\Lambda$ there exists a full embedding of the category
of finite dimensional $\Lambda-$modules into ${\cal VB}(X)$. \\ 
Let $X$ be an algebraic curve. How to normalize it? There are
several methods, algorithms and implementations for this 
purpose. A new algorithm and implementation is presented 
in~\cite{DJGP}.

\subsection{ Foliations}~
A simplest example of foliation is a trivial $k$ dimensional
foliation or a trivial codimension $n-k$ foliation of Euclidean
space ${\bf R}^{n}={\bf R}^{k} \times {\bf R}^{n-k}, 0 \le k \le
n:$
$$ {\bf R}^{n} = \bigcup_{(x_{k+1},\ldots ,x_{n})}{\bf R}^{k}
\times (x_{k+1},\ldots,x_{n}), $$
that is, $ {\bf R}^{n}$ decomposes into a union of
$ {\bf R}^{k} \times (x_{k+1}, \ldots ,x_{n}) $ 's each of
which is $ C^{\infty}$ diffeomorphic to $ {\bf R}^{k}.$
$ {\bf R}^{k} \times (x_{k+1}, \ldots ,x_{n}) $ is called
a leaf of the foliation. \\
Still one example can be obtained from consideration of
nonsingular vector fields on the torus~\cite{Ta,CS}. More generally,
a nonsingular flow on a manifold  corresponds to a foliation
on the manifold by one dimensional leaves where the leaves
are provided with Riemannian metrics and directed. \\
Consider vector field
\begin{equation}
 \label{VF}
      Q(x,y)\partial x + P(x,y)\partial y, 
\end{equation}         
where $ P $ and $ Q $ are complex polynomials of degree $ n $ in two
variables. The vector field (or corresponding dual 1-form
${\omega} = P(x,y)dx -Q(x,y)dy $) gives rise to a foliation ${\cal F} 
$ of degree $ n $ of two dimensional projective space 
$ {\bf P}^{2} $ over $ {\bf C} $
by Riemann surfaces and singular points. It is naturally to ask
about limit cycles and multivalued first integrals of (\ref{VF})  
and a topology of the foliation.
   The investigation of (\ref{VF}) began in~\cite{Pet}, where
the case $ n = 2 $ considered. For  the investigation of (\ref{VF}) 
and ${\cal F} $ it is naturally to
introduce of algebraic geometric methods and algebraic geometric
invariants~\cite{Shap,Bif}. These include blow-ups, divisors,
indexes of singular points, Chern classes of vector bundles over some
Riemann surfaces and computation of holonomy (monodromy) groups.
Already in the case of~\cite{Pet} the complex dimension of the
space of coefficients of (1)
is equal 12 (real dimension is equal 24). This defines the
application of computer algebra in the case. \\

\subsubsection{Kronecker foliation }~\cite{Co:NG}
  Let ${\theta}$ be a irrational number. Consider a vector 
field $ F = \partial x + \theta\partial y $ on the two
dimensional torus 
$ T^2 = {\bf R}/{\bf Z} \times {\bf R}/{\bf Z}.$
The Kronecker foliation ${\cal F}_{\theta}$ is defined
by the flow of $F$ in the following way: each leaf is
labeled by a point of circle $S^{1} = {\bf R}/{\bf Z}.$ 
So the leaf space is the foliated bundle
$$ ({\bf R} \times {\bf R}/{\bf Z})/{\bf Z}, \;
(x,y) \sim (x + 1,y + \theta),$$
where the equivalence relation on $(x,y) \in {\bf R} \times S^{1}$
is defined by
$$ R_{\theta}(x,y) = (x,y + \theta) \bmod 1 .$$
 By the foliation ${\cal F}_{\theta}$
the $C^*$-algebra ${\cal A}_{\theta}$ is constracted.
${\cal A}_{\theta}$ is generated by a pair of unitary
symbols subject to the relation
\begin{equation}
 \label{TA}
  UV = \exp (2\pi i \theta) VU.
\end{equation}

\subsubsection{ Formal Groups}
  Formal groups was introduced by J. Dieudonne and M. Lazard 
around 1954. During 1968-71 and later were found interesting 
connections between topology and formal groups~\cite{Qu,BMN}. 
Also was found connection
between formal groups and zeta functions~\cite{Ho}. 
There is a connection of formal groups with characteristic
classes of foliations and $K-$theory~\cite{No}. 
The construction of the Grothendiesk group $K_{0}$ for algebraic 
monoid is rather simple~\cite{Ka:K}. What is the $K_{0}$ functor 
in the case of formal groups? Let us consider a simple example. \\
{\bf Proposition.} Let ${\cal FG}_k$ be the category of commutative 
formal groups over field $k$ such that every formal group in the 
category is the product of one dimensional formal groups. Then 
$K_{0}({\cal FG}_{k})$ is a free abelian group with infinite 
number of generators. \\ 
                               
\subsubsection{Moduli spaces }
   The theory of moduli spaces~\cite{Ma:T,HM:MC} has, in recent
years, become the meeting ground of several different branches
of mathematics and physics - algebraic geometry, instantons, 
differential geometry, string theory and arithmetics. Here 
we recall some underlieing algebraic structures of the relation.
In previous subsection we have reminded the situation with vector
bundles on projective algebraic curves $X$. On $X$ any first Chern
class $c_{1} \in H^{2}(X,{\bf Z})$ can be realized as $c_1$ of 
vector bundle of prescribed rank (dimension) $r.$  How to classify
vector bundles over algebraic varieties of dimension more than $1?$
This is one of important problems of algebraic geometry and the
problem has closed connections with gauge theory in physics and 
differential geometry. Mamford~\cite{Ma:T} and others have 
formulated the problem about the determination of which cohomology
classes on a projective variety can be realized as Chern classes
of vector bundles?  Moduli spaces are appeared in the problem.
What is moduli? Classically Riemann claimed that $3g - 3$ (complex) 
parameters could be for Riemann surface of genus $g$ which would
determine its conformal structure (for elliptic curves, when
$g = 1,$ it is needs one parameter). From algebraic point of view 
we have the following problem: given some kind of variety,
classify the set of all varieties having something in common with
the given one (same numerical invariants of some kind, belonging
to a common algebraic family). For instance, for an elliptic
curve the invariant is the modular invariant of the elliptic
curve.  \\
Let {\bf B} be a class of objects. Let $ S $  be a scheme. 
A family of objects parametrized by
the $ S $ is the  set of objects

    $$ X_{s}: s \in S, X_{s} \in {\bf B} $$

 equipped with an additional structure compatible with the structure
  of the base $ S $.
Parameter varieties is a class of moduli spaces. These varieties
is very convenient tool for computer algebra investigation
of objects that parametrized by the parameter varieties. We
have used the approach for investigation of rational points
of hyperelliptic curves over prime finite fields~\cite{Gl:MS}.

\subsection {Connection      }
  Consider the connection in the context of algebraic geometry.
We shall base at one I. Shafarevitch's seminar on algebraic geometry.
Let $S/k$ be the smooth scheme over field $k$, $U$ an element of
open covering of $S$, $ {\cal O}_S $ the structure sheaf on $S$,
$\Gamma(U,{\cal O}_{S})$ the sections of ${\cal O}_{S}$ on $U$. 
Let $\Omega^{1}_{S/k}$ be the sheaf of germs of $1-$dimension 
differentials, $\cal F$ be a coherent sheaf. The {\em connection} on the sheaf 
$\cal F$ is the sheaf homomorphism 
$$ \nabla: {\cal F} \rightarrow \Omega^{1}_{S/k} \otimes {\cal F},$$
such that, if $ f \in \Gamma(U,{\cal O}_{S}), \;
g \in \Gamma(U,{\cal F})$ then
$$ \nabla(fg) = f\nabla(g) + df \otimes g.$$ 
  There is the dual definition. Let ${\cal F}$ be the locally
free sheaf, $\Theta^{1}_{S/k}$ the dual to sheaf $\Omega^{1}_{S/k}$, 
$\partial \in \Gamma(U,\Theta^{1}_{S/k})$. 
The {\em connection} is the homomorphism  
$$ \rho: \Theta^{1}_{S/k} \rightarrow 
End_{{\cal O}_S}({\cal F},{\cal F}). \;
\rho(\partial)(fg) = \partial(f)g + f\rho(\partial). $$ 

\subsubsection {Integration of connection} 

Let $\Omega^{i}_{S/k}$ be the sheaf of germs of $i-$differentials,
$$ \nabla^{i}(\alpha \otimes f) = d\alpha \otimes f + (-1)^{i}\alpha
\wedge \nabla(f) $$
   Then $ \nabla, \; \nabla^{i}$ define the sequence of homomorphisms:
\begin{equation}
 \label{CC}
 {\cal F} \rightarrow \Omega^{1}_{S/k} \otimes {\cal F}
 \rightarrow \Omega^{2}_{S/k} \otimes {\cal F} \rightarrow \cdots     , 
\end{equation}  
  A connection is {\em integrable} if (\ref{CC}) is a complex.  
This is equivalent to $\nabla \otimes \nabla^{1} = 0.$

\section{ Zeta Functions }

  We have considered some validated numerics aspects of evaluation 
of zeta functions in papers~\cite{G:VN,G:IE}. Here we shall consider 
some computer algebra aspects of computation and evaluation of values
of zeta functions.
\subsection{Riemann zeta}

 Consider the series  $\zeta(s) = \sum_{n} \frac{1}{n^s}$
for complex values of $s$ with $Re(s) \geq 1.$ Analytical
formula for computation of values of the $\zeta(s)$ can be 
taken from the result of Backlund (a little bit reformulated). 
The result is based of Euler-Maclaurin summation.\\
{\bf Proposition} (Backlund) Let $N$ be natural $> 1.$
Let $s = \sigma + it$ and let $\sigma \geq 1.$ Let
$B_{2k}$ be the Bernoulli numbers in even numeration, \\
$ S(N-1,s) = \sum_{n=1}^{N-1} n^{-s} + \frac{N^{1-s}}{s-1},
\\
 B(N,k,s)
  = \frac{1}{2}N^{-s} + \frac{B_2}{2}sN^{-s-1} + \cdots
+ \frac{B_{2k}}{(2k)!}s(s+1) \ldots (s+2k-2)N^{-N-2k+1} .$ \\
Then
\begin{equation}
 \label{EMS}
\zeta(s) = S(N-1,s) + B(N,k,s) + R_{2k},
\end{equation}
where $$|R_{2k-2}| \leq \left| \frac{s+2k-1}{\sigma + 2k -1}
\right| |B_{2k} \; \mbox{term of }(\ref{EMS})|. $$

\subsection{Dirichlet series}

Let $K = {\bf Q}(\sqrt{D})$ be a real quadratic field with
positive integer squarefree $D$ and $\chi(n) =
\left( \frac{\Delta}{n} \right)$ be the Kronecker symbol. Here
$$\Delta = \left\{  \begin{array}{ll}
                        D, & \mbox{ $ D \equiv 1 \pmod 4 $ ,} \\
                        4D, & \mbox{ $ D \equiv 2,3 \pmod 4  $.}
                        \end{array}
                \right. $$
Let $A = \pi/\Delta, \; E(x) = \int_{x}^{\infty} \frac{e^{-t}}{t}dt,
\; erfc(x) = \frac{2}{\pi} \int_{x}^{\infty} e^{-t^2}dt. $
  \\
{\bf Proposition}~\cite{WB,PB}. \\
$$ L(1,\chi) = \frac{1}{\sqrt{\Delta}} \sum_{n=1}^{m} \chi(n) E(An^2)
\; + \sum_{n=1}^{m} \left( \frac{\chi(n)}{n} \right)
erfc(n \sqrt{A}) + R_m,$$
where $|R_m| < \frac{\Delta^{3/2}}{\pi^2} \frac{e^{-Am^2}}{m^3}.$

\subsection{Remarks about $L-$functions of elliptic curves}~\cite{Ta:AE}.
 Let $E/{\bf Q}$ be an elliptic curve given in Weierstrass form
by an equation
\begin{equation}
\label{WN}
E: y^2+a_1 xy + a_3 y = x^3 + a_2 x^2 + a_4 x + a_6 ,
\end{equation}
and let $b_{2} , b_4, b_6, b_8, c_4, c_6, \Delta , j$ be the
usual associated quantities~\cite{Ta:AE}. Let now~(\ref{WN}) be a
global minimal Weierstrass equation for $E$ over {\bf Z.} For each
prime $p$ the reduction~(\ref{WN}) $\pmod p$ defines a curve $E_p$ over
the prime field ${\bf F}_p$. Let $A_p$ denote the number of
 points of $E_p$ rational over ${\bf F}_p$. Let
$$ t_p = 1 + p - A_p. $$
If $p \not| \Delta,$ then $t_p$ is the trace of Frobenius
and satisfies $|T_p| \leq 2\sqrt{p}.$ In the case
{\it Artin-Hasse zeta function} of the elliptic curve $E_p$
is:
\begin{equation}
\label{ArH}
\zeta_{E_p}(s) =
\frac{1 - t_{p}p^{-s} + p^{1-2s}}{(1 - p^{-s})(1 - p^{1 - s})}.
\end{equation}
 If $p \mid \Delta,$ then $E_p$ is not an elliptic curve and has
a singularity $S.$  In the case
$$ t_p = \left\{ \begin{array}{ll}
                        0, & \mbox{if $S $ is a cusp,} \\
                        1, & \mbox{if $S $ is a node,} \\
                       -1, & \mbox{if $S $ is a node with tangent quadratic
                over ${\bf F}_p$}
                        \end{array}
                \right. $$
The {\it Hasse-Weil } $L-$function of $E/{\bf Q}$ is defined by equation
\begin{equation}
\label{HW}
\L_{E}(s) = \prod_{p \, | \Delta}
\frac{1 }{(1 - t_p p^{-s})} \prod_{p \, \not | \, \Delta}
\frac{1}{1 - t_{p}p^{-s} + p^{1-2s}}.
\end{equation}
   From the work of A. Wiles, R. Taylor and A. Wiles, and work
of F. Diamond it is known that (semistable) elliptic curves
over $E/{\bf Q}$  are modular. Knowing the modularity of
$E/{\bf Q}$  is equivalent to the existence of a modular form $f$
on $\Gamma_{0}(N)$ for some natural value $N,$ which we write
$f = \sum a_n q^n.$ The $L-$function of $E$ is thus given by
the Mellin transform of $f, \; L(f,s) = \sum a_n q^n.$ In 
particular, the behavior of $L(E,s)$ at $s = 1$ can be deduced
from modular properties of $E.$ \\
Let $E/{\bf Q}$ be a modular elliptic curve and the global minimal
model of the $E/{\bf Q}$ has prime conductor $l$. Let $p$ be a prime
and $A_p$  be the number of points of $E_p$ in 
${\bf F}_p.$ Then there is exists a modular form $f$ on
$\Gamma_{0}(l), \;  f = \sum_{n=1}^{\infty} a_n q^n,$ where
$a_p, \; p \neq l,$ equals $p + 1 - A_p.$

\section{Some Algebraic Methods and Structures of 
Noncommutative  Geometry }
At first remind two definitions. \\
{\em Hopf algebras}. \\
  Let $\cal H$ be an algebra with unit $e$ over field $k.$
Let $a, \; b \in {\cal H}, \; \mu(a,b) = ab$ the product in 
$\cal H$, $\alpha \in k, \; p:k \rightarrow {\cal H},
\; p(\alpha) = \alpha e $ the unit, $\varepsilon: {\cal H} 
\rightarrow k, \; \varepsilon(a) = 1, $ 
the counit, $\Delta(a) = a \otimes a$ the coproduct, $S: \; {\cal H} 
\rightarrow {\cal H}$  the antipode, such that the axioms \\
1)$(1 \otimes \Delta) \circ \Delta = (\Delta \otimes 1) \circ \Delta$ 
(coassociativity); \\
2)$\mu \circ (1 \otimes S) \circ \Delta = \mu \circ 
(S \otimes 1) \circ \Delta = p \circ \varepsilon$ (antipode). \\ 
Then the system $({\cal H}, \mu, p, \varepsilon, \Delta, S)$
is called the {\em Hopf algebra}.        \\
 Let ${\cal A}$ and ${\cal B}$ be
$C^*$-algebras. Let ${\cal K}$ be the $C^*$-algebra
of compact operators. If ${\cal A} \otimes {\cal K}$ is
isomorphic to ${\cal B} \otimes {\cal K}$ then ${\cal A}$
is called {\em Morita equivalent} to ${\cal B}$. \\
  We give here some remarks to algebraic setting of
papers~\cite{Co:N2,Co:NG,Mo:OA}. In a general context this algebraic
setting includes: \\
-\, Operator algebras. \\
-\, Representation theory. \\
-\, K-theory. \\
-\, Algebraic geometry. \\
  More specificaly: \, $C^*$-algebras; finite projective modules; 
holonomy groupoid; groupoid for transformation groups;
group $C^*$-algebras; Morita equivalence; K-theory; von Neumann 
algebras; cyclic cohomology; invariant transversal measures and the
Ruelle-Sullivan current; Godbillon-Vey class; Hopf algebras. \\
   One of the paradigm of noncommutative geometry is to describe
the geometry of ordinary space in terms of the algebra of functions
and then deforming to the noncommutative case. By deformation theory
as the part of moduli space theory some constructions of section 1
can be embedded into noncommutative geometry.

\section{Appendix. The Common Lisp text of the packages CUTSET and
CUTSETDG} 
  Rooted graph $ G $ is the graph that each node of $ G $
is reachable from a node $ r \in G .$
By {\it cutset} of the graph we shall understand an
appropriated subset of nodes (called cutpoints) such
that any cycle of the graph contains at least one
cutpoint. DFS is the abbreviation of Depth First Search
method. During DFS we numbering nodes and label (mark)
edges. By Tarjan~\cite{Ta} the DFS method has linear complexity.
The packages are implemented on Allegro CL 3.0.2~\cite{G,CL}. 
  The main function Cutsetdg of the package CUTSETDG
computes Cutset (a subset of
vertices which cut all cycles in the graph) of arbitrary
rooted directed graph. It uses 3 functions: Adjarcn,
Cutpoints and Unicut. For this program the author developed
rather simple and efficient algorithm that based on
DFS-method. The function is the base for Cutset methods
for systems of procedures.
    The main function Cutset of the Package CUTSET
computes Cutset of arbitrary
rooted graph. It also uses 3 functions: Adjedgn,
Cutpoints and Unicut. An algorithm for this function that
based on DFS-method is also developed by the author. The algorithm
is rather simple and efficient. The function is the base
for Cutset methods.
Package CUTSET    \\
 Title: Cutset of rooted graph  \\
 Summery:
 This package implements computation of a cutset
 of rooted graph. The graph have to defined by adjacency
 list. \\
 Example of Call:   (Cutset 'a)
 where "a" is a root of exploring graph. \\

\begin{verbatim}
;; Name: Cutset
;;
;; Title: Cutset of rooted graph
;;
;; Author: Nikolaj M. Glazunov
;;
;; Summery:
;; This package implements computation of a cutset
;; of rooted graph. The graph have to defined by adjacency
;; list
;;
;; Allegro CL 3.0.2
;;
(DEFUN Cutset (startnode)

;; startnode is a root of exploring graph
;; DFS - Depth First Search method for connected graph
;; cuts - cutset of the graph
;; st - stack
;; v - exploring node
;; sv -  son of the exploring node
;; inl - inverse edges list
;; df -  DFS-numbering of the node (property)
;; Adjed - adjacency list of the node (property)
;; 1, 2  - lables

   (SETF (GET startnode 'df ) 1)   ;DFS-numbering is equal 1
         (SETF (GET startnode 'Adjed)
            (Adjedgn  (LIST startnode))) ;;startnode obtained the a-list
                                 ;; property Adjed (adjacency edges)

   (PROG ( st v sv inl  cuts)
         (PUSH startnode st)
   1   (SETQ v (CAR st)) ;;explored node received value from stack of nodes

   2      (COND ((NOT (EQ  NIL (GET v 'Adjed)))
                    (SETQ sv (CAAR (GET v 'Adjed)))

             (COND ((EQ NIL (GET sv 'df))
                          (SETF (GET sv 'df) (+ 1 (GET v 'df)))
       ;; modification of the edge of son's adj-list
                      (SETF (GET sv 'Adjed)
            (Adjedgn  (LIST sv))) ;; node sv obtained the adj-list property
                                 ;; Adjed (adjesency edges)
     (SETF (GET v 'Adjed) (REMOVE (LIST sv v) (GET v 'Adjed) :TEST 'EQUAL))
     (SETF (GET sv 'Adjed) (REMOVE (LIST v sv) (GET sv 'Adjed) :TEST 'EQUAL))
                     (PUSH sv st)
                     (GO 1)
            )
               (T   ;; the node has number.
                       ;; the edge is inverse
                         (SETQ inl (CONS (CAR (GET v 'Adjed)) inl))
     (SETF (GET v 'Adjed) (REMOVE (LIST sv v) (GET v 'Adjed) :TEST 'EQUAL))
     (SETF (GET sv 'Adjed) (REMOVE (LIST v sv) (GET sv 'Adjed) :TEST 'EQUAL))
                     (GO 2)
                    )
               )
                 )
                           (T
                (POP st)
                  (COND     ((NOT (NULL st))
                      (GO 1)
              )
            (T                 ;; end
                   (SETQ cuts (Unicut (Cutpoints inl)))
                        (RETURN cuts)
                )
              )
                ))
     )
   )
(DEFUN Adjedgn (PATH)
      (MAPCAN #'(LAMBDA (E)
              (COND ((MEMBER (CAR E) (CDR E)) NIL)
                   (T (LIST E))))
          (MAPCAR #'(LAMBDA (E)
                            (CONS E PATH))
                (GET (CAR PATH) 'NEIGHBORS))))

(DEFUN Cutpoints (inl)
        (MAPCAR 'CAR inl))
(DEFUN Unicut (cpl)
     (COND ((NULL cpl) NIL)
           (T (CONS (CAR cpl)
                (Unicut (REMOVE (CAR cpl) cpl)))
             )
     )
     )
\end{verbatim}

Package CUTSETDG \\
 Title: Cutset of a rooted directed graph   \\
 Summery:
 This package implements computation of a cutset
 of a rooted directed graph. The graph have to defined by adjacency
 list. \\
 Example of Call:   (Cutsetdg 'a)
 where "a" is a root of exploring graph. \\

{\bf The Common Lisp text of the package CUTSETDG} \\
\begin{verbatim}

;; Name: Cutsetdg
;;
;; Title: Cutset of a rooted directed graph
;;
;; Author: Nikolaj M. Glazunov
;;
;; Summery:
;; This package implements computation of a cutset
;; of a rooted directed graph. The graph have to defined by adjacency
;; list
;;
;; Example of Call:   (Cutsetdg 'a)
;; where "a" is a root of exploring graph
;;
;; Allegro CL 3.0.2
;;
(DEFUN Cutsetdg (startnode)

;; startnode is a root of exploring graph
;; DFS - Depth First Search method for connected graph
;; st - stack
;; v - exploring node
;; sv -  son of the exploring node
;; inl - inverse edges list
;; df -  DFS-numbering of the node (property)
;; Outarcs - adjacency list of Output arcs of the node (property)
;; 1, 2  - lables

    (SETF (GET startnode 'df ) 1)   ;DFS-numbering is equal 1
   (SETF (GET startnode 'Outarcs)
            (Adjarcn  (LIST startnode))) ;;startnode obtained the a-list
                                 ;; property Outarcn (Output arcs of the
                                 ;; node)
    (PROG ( st v sv inl  cuts)
         (PUSH startnode st)
   1   (SETQ v (CAR st)) ;;explored node received value from stack of nodes
   2      (COND ((NOT (EQ  NIL (GET v 'Outarcs)))
                    (SETQ sv (CAAR (GET v 'Outarcs)))
                        (COND ((EQ NIL (GET sv 'df))
                            (SETF (GET sv 'df) (+ 1 (GET v 'df)))
       ;; modification of the edge of son's adj-list
       (SETF (GET sv 'Outarcs)
            (Adjarcn  (LIST sv))) ;; node sv obtained the adj-list property
                                 ;; Outarcs (adjacency Outarcs)
          (SETF (GET v 'Outarcs) (REMOVE (LIST sv v) (GET v 'Outarcs)
                              :TEST 'EQUAL))
                           (PUSH sv st)
                            (GO 1)
            )
              (T
                 (COND ((AND (> (GET v 'df) (GET sv 'df))
                            (NOT (EQUAL v sv)))
              (SETQ inl (CONS (CAR (GET v 'Outarcs)) inl))
              (SETF (GET v 'Outarcs) (REMOVE (LIST sv v) (GET v 'Outarcs)
                              :TEST 'EQUAL))
                        (GO 2)
                        )
              (T
               (SETF (GET v 'Outarcs) (REMOVE (LIST sv v) (GET v 'Outarcs)
                              :TEST 'EQUAL))
                (GO 2)
                )
                       )
                )
                   )
                 )
                              (T
                (POP st)
         (COND     ((NOT (NULL st))
                      (GO 1)
              )
                     (T                 ;; end
                   (SETQ cuts (Unicut (Cutpoints inl)))
                (RETURN cuts)
                )
              )
                ))
     )
   )
(DEFUN Adjarcn (PATH)
      (MAPCAN #'(LAMBDA (E)
              (COND ((MEMBER (CAR E) (CDR E)) NIL)
                   (T (LIST E))))
          (MAPCAR #'(LAMBDA (E)
                            (CONS E PATH))
             (CAR  (GET (CAR PATH) 'NEIGHBORS)))))
(DEFUN Cutpoints (inl)
        (MAPCAR 'CAR inl))
(DEFUN Unicut (cpl)
     (COND ((NULL cpl) NIL)
           (T (CONS (CAR cpl)
                (Unicut (REMOVE (CAR cpl) cpl)))
             )
     )
     )
\end{verbatim}

\end{document}